\documentclass[reqno]{amsart}
\usepackage{amsrefs}
\usepackage{amssymb}

\newtheorem{thm}{Theorem}
\newtheorem{cor}[thm]{Corollary}

\newcommand{\abs}[1]{\left\lvert#1\right\rvert}
\newcommand{\solid}[0]{\to}
\newcommand{\flimsy}[0]{\leadsto}
\newcommand\fixspacetop{\rule{0pt}{2.2ex}}
\DeclareMathOperator{\ord}{ord}

\begin{document}

\title{Wieferich pairs and Barker sequences, II}
\author{Peter Borwein}
\address{Department of Mathematics\\
Simon Fraser University\\
Burnaby, BC V5A~1S6 Canada}
\thanks{Research of P.~Borwein supported in part by NSERC of Canada and MITACS}
\email{pborwein@sfu.ca}
\author{Michael J. Mossinghoff}
\address{Department of Mathematics, Davidson College, Davidson,
NC 28035-6996 USA}
\thanks{Research of M.~J. Mossinghoff was partially supported by a grant from the Simons Foundation (\#210069).}
\email{mimossinghoff@davidson.edu}

\date\today
\subjclass[2000]{Primary: 05B20, 94A55; Secondary: 05-04, 05B10, 11A07}
\keywords{Barker sequence, circulant Hadamard matrix, Wieferich prime pair.}

\begin{abstract}
We show that if a Barker sequence of length $n>13$ exists, then either $n=3\,979\,201\,339\,721\,749\,133\,016\,171\,583\,224\,100$, or $n > 4\cdot10^{33}$.
This improves the lower bound on the length of a long Barker sequence by a factor of nearly $2000$.
We also obtain $18$ additional integers $n<10^{50}$ that cannot be ruled out as the length of a Barker sequence, and find more than 237000 additional candidates $n<10^{100}$.
These results are obtained by completing extensive searches for Wieferich prime pairs and using them, together with a number of arithmetic restrictions on $n$, to construct qualifying integers below a given bound.
We also report on some updated computations regarding open cases of the circulant Hadamard matrix problem.
\end{abstract}

\maketitle

\section{Introduction}\label{sectionIntroduction}

The \textit{$k$th aperiodic autocorrelation} for a finite sequence $a_1$, \ldots, $a_n$ is defined by
\[
c_k = \sum_{i=1}^{n-k} a_i a_{i+k}.
\]
A \textit{Barker sequence} of length $n$ is a sequence $a_1$, \ldots, $a_n$, each $\pm1$, with the property that each of its aperiodic autocorrelations besides $c_0$ is small: one requires that $\abs{c_k}\leq1$ for each $k\geq1$.
Barker sequences are known only for lengths $n\in\{1,2,3,4,5,7,11,13\}$, and in fact there is only one Barker sequence of each of these lengths, after accounting for certain symmetries.
(Reversing, negating, and negating every other term of a Barker sequence always produces another Barker sequence.)
It is widely conjectured that no other Barker sequences exist.
Turyn and Storer \cite{TS} proved that $n=13$ is the maximal length of a Barker sequence of odd length, and many restrictions are known for the even case.

The even case of the Barker sequence problem is related to another well-known problem in combinatorial optimization.
For a sequence $a_1$, \ldots, $a_n$, we first define its \textit{$k$th periodic autocorrelation} by
\[
\gamma_k = \sum_{i=1}^{n} a_i a_{i+k},
\]
where the indices on the right are taken modulo $n$, so that $\gamma_k = c_k + c_{n-k}$ for each $k\geq0$.
It is well known \cite{EK92} that if $a_1$, \ldots, $a_n$ is a Barker sequence with $n>2$, then its sequence of periodic autocorrelations is constant for $k\geq1$, and that their common value $\gamma$ is $0$ if $n$ is even, $1$ if $n\equiv1$ mod $4$, and $-1$ if $n\equiv3$ mod $4$.
It follows that if $a_1$, \ldots, $a_n$ is a Barker sequence and $n>13$, then the vector $\mathbf{a}=(a_1,\ldots,a_n)$ is orthogonal to all nontrivial cyclic shifts of itself.
If we record each successive cyclic shift of $\mathbf{a}$ as a row of an $n\times n$ matrix $H$, then $H$ is a circulant matrix with mutually orthogonal rows in which each entry is $\pm1$.
Such a matrix is a \textit{circulant Hadamard matrix}, and it is widely conjectured that no such matrix exists with order $n>4$.

A number of restrictions on permissible values of $n$ are known in the circulant Hadamard matrix problem.
We review these in Section~\ref{secRestrictions}, along with an additional arithmetic restriction on $n$ in the open case of the Barker sequence problem.
In 2009, the second author \cite{Moss09} employed these restrictions to determine the smallest integer $n>13$ that could not be eliminated as the length of a Barker sequence, proving that either
\begin{equation}\label{eqnN0}
n = 189\,260\,468\,001\,034\,441\,522\,766\,781\,604
\end{equation}
or $n>2\cdot10^{30}$.
The same article established that fewer than $1600$ integers $n<4\cdot10^{26}$ satisfy all the known restrictions in the circulant Hadamard matrix problem.
Recently, Leung and Schmidt \cite{LS12} obtained some further restrictions on the order of a circulant Hadamard matrix, and some of these restrictions apply to the Barker sequence problem as well.
These new restrictions are also summarized in Section~\ref{secRestrictions}.
Using these new conditions, Leung and Schmidt proved that the exceptional value \eqref{eqnN0} is not a permissible length for a Barker sequence.

In this article, we determine the smallest integer $n$ that satisfies all of the known conditions now required for $n>13$ to be the length of a Barker sequence, and we show in Section~\ref{secComputations} that this is the only such integer less than $4\cdot10^{33}$.
Our method follows that of \cite{Moss09}, but incorporates the new restrictions of \cite{LS12}.
This method relies on a large search for \textit{Wieferich prime pairs}, which are pairs of prime numbers $(q,p)$ with the property that $q^{p-1}\equiv1$ mod $p^2$.
In Section~\ref{secAdditional}, we determine more than $500$ additional integers that satisfy all of the known requirements for the length of a Barker sequence, including eighteen more values less than $10^{50}$.
We also compile more than $237000$ integers less than $10^{100}$ which survive most of the known restrictions, but for which complete testing is at present computationally prohibitive.
We expect, however, that the vast majority of these integers will in fact satisfy all of the known requirements.

In particular, we prove the following theorem.

\begin{thm}\label{thmBarkerBound}
If $n>13$ is the length of a Barker sequence, then either
\[
n = 3\,979\,201\,339\,721\,749\,133\,016\,171\,583\,224\,100,
\]
or $n>4\cdot10^{33}$.
\end{thm}

Here, $n=4u^2$ with $u=5 \cdot 13 \cdot 29 \cdot 41 \cdot 2953 \cdot 138200401$; the prior value \eqref{eqnN0} is $4u_0^2$ with $u_0=13 \cdot 41 \cdot 2953 \cdot 138200401$.

Last, in Section~\ref{secCHM} we update some computations on the circulant Hadamard matrix problem, and show that the number of open cases with $n\leq4\cdot10^{26}$ may be reduced from $1393$ to $1371$.

\section{Arithmetic restrictions}\label{secRestrictions}

We review the known arithmetic restrictions on the order of a circulant Hadamard matrix, since these are automatically restrictions for the Barker sequence problem in the open case where the length is even.
We also review a requirement specifically for Barker sequences.

First, it is well known that if $n$ is the order of a Hadamard matrix and $n>2$, then $4\mid n$, and that the order of a circulant Hadamard matrix is a square, so we write $n=4u^2$.
Turyn \cite{Turyn65} proved that $u$ must be odd, and cannot be a prime power.
He established an important general criterion as well, known as the self-conjugacy test.
We require two brief definitions for this.
First, for integers $r$ and $s$, we say that $r$ is \textit{semiprimitive} mod $s$ if there exists an integer $j$ such that $r^j\equiv -1$ mod $s$.
Necessarily, then, $\gcd(r,s)=1$.
Second, we say that $r$ is \textit{self-conjugate} mod $s$ if every prime divisor $p\mid r$ is semiprimitive modulo the $p$-free part of $s$, that is, the largest divisor of $s$ that is not divisible by $p$.
We may now state Turyn's self-conjugacy test.

\begin{thm}\label{thmTuryn}
Suppose $n=4u^2$ is the order of a circulant Hadamard matrix, and $r$ and $s$ are integers with $r\mid u$, $s\mid n$, and $\gcd(r,s)$ has $k\geq1$ distinct prime factors.
If $r$ is self-conjugate mod $s$, then $rs\leq2^{k-1}n$.
\end{thm}

By taking $r=p^\ell$ and $s=2p^{2\ell}$ in the self-conjugacy test, one obtains a special case, noted for instance in \cite{JL}.
It states that the order of a circulant Hadamard matrix cannot have a prime-power factor that is too large.

\begin{cor}\label{corLargePrime}
If $n=4u^2$ is the order of a circulant Hadamard matrix, and $p^\ell\mid u$ for an odd prime $p$ and positive integer $\ell$, then $p^{3\ell}\leq2u^2$.
\end{cor}

In \cite{LS12}, Leung and Schmidt obtained another restriction that rules out values of $n$ having a sizable prime-power divisor, provided that a side condition also holds.
For a prime $p$ and integer $t$, let $\nu_p(t)$ denote the largest integer $k$ such that $p^k\mid t$.

\begin{thm}\label{thmLS5}
Suppose that $n=4u^2$ is the order of a circulant Hadamard matrix, let $p$ be an odd prime dividing $u$, let $a=\nu_p(u)$, and suppose that $p^{2a}>2u$.
Further, let $r$ be a divisor of $m=u/p^a$, with $r$ self-conjugate modulo $p$, and suppose that $q_1$, \ldots, $q_k$ are the prime divisors of $m/r$.
Then
\[
\gcd(\ord_p(q_1),\ldots,\ord_p(q_k)) \leq \frac{m^2}{r^2}.
\]
\end{thm}

This result rules out \eqref{eqnN0} as the possible length of a Barker sequence, by taking $p=138200401$ and $r=2953$.

Next are two tests based on the field descent method of Schmidt \cite{Schmidt99}.
For this, let $\mathcal{D}(t)$ denote the set of prime divisors of the integer $t$.
For a positive integer $m$ and a prime $q$, let
\[
m_q = \begin{cases}
\displaystyle\prod_{p\in\mathcal{D}(m)\backslash\{q\}} p, & \textrm{if $m$ is odd or $q=2$},\\
\displaystyle2\prod_{p\in\mathcal{D}(m)\backslash\{q\}} p, & \textrm{otherwise}.
\end{cases}
\]
In the first case, then, $m_q$ is simply the $q$-free and squarefree part of $m$.
Next, let $\ord_s(t)$ denote the order of $t$ in the multiplicative group of units modulo $s$, and for positive integers $m$ and $n$ and a prime $r$, define $b(r,m,n)$ by
\[
b(r,m,n) = \begin{cases}
\displaystyle\max_{q\in\mathcal{D}(n)\backslash\{2\}} \left\{\nu_2(q^2-1) + \nu_2(\ord_{m_q}(q)) - 1\right\}, & \textrm{if $r=2$},\\
\displaystyle\max_{q\in\mathcal{D}(n)\backslash\{r\}} \left\{\nu_r(q^{r-1}-1) + \nu_r(\ord_{m_q}(q))\right\}, & \textrm{if $r>2$},
\end{cases}
\]
with the convention that $b(2,m,2^k)=2$ and $b(r,m,r^k)=1$ for an odd prime $r$.
Finally, define $F(m,n)$ by
\[
F(m,n) = \gcd\left(m, \prod_{p\in\mathcal{D}(m)} p^{b(p,m,n)}\right).
\]
From \cite{LS05}*{Cor.\ 4.5}, we have the following highly restrictive inequality.

\begin{thm}\label{thmLS}
If $n=4u^2$ is the order of a circulant Hadamard matrix, then $u\varphi(u) \leq F(u^2,u)$.
\end{thm}

In \cite{LS12}, Leung and Schmidt established a second bound depending on the same function $F(m,n)$.
We cite their result here only as it applies to circulant Hadamard matrices.

\begin{thm}\label{thmLS1}
If $n=4u^2$ is the order of a circulant Hadamard matrix, and $m$ and $w$ are positive integers with $m\mid u$, $w\mid n$, and $m$ is self-conjugate modulo $n/w$, then $n\varphi(F(n/w,u^2/m^2))) \leq w^2 F(n/w,u^2/m^2)^2$.
\end{thm}

Finally, Eliahou, Kervaire, and Saffari \cite{EKS} established the following restriction specifically for Barker sequences.

\begin{thm}\label{thmEKS}
If $n=4u^2$ is the length of a Barker sequence, and $p$ is a prime number with $p\mid u$, then $p\equiv1$ mod $4$.
\end{thm}

We remark that Leung and Schmidt in \cite{LS12} obtain a third new restriction for circulant Hadamard matrices of order $n=4u^2$, for the case where every prime divisor $p$ of $u$ satisfies $p\equiv3$ mod $4$.
Since this restriction is not relevant to the Barker sequence problem, we do not consider it here, and instead discuss it in Section~\ref{secCHM}.

\section{Proof of Theorem~\ref{thmBarkerBound}}\label{secComputations}

We employ the method of \cite{Moss09} to determine all integers $u\leq U = 10^{16.5}$ composed entirely of primes congruent to $1$ mod $4$ and less than $P=2^{1/3}\cdot10^{11}$, and for which $u\varphi(u)\leq F(u^2,u)$.
This allows us to construct all integers $n=4u^2$ up to $4\cdot10^{33}$ which satisfy Corollary~\ref{corLargePrime}, Theorem~\ref{thmLS}, and Theorem~\ref{thmEKS}.
Note that in this case
\[
F(u^2,u)=\gcd\bigl(u^2,\prod_{p\mid u} p^{b(p,u^2,u)}\bigr)
\]
and, since $u$ must be odd and we can assume that $u$ is not a prime power, that
\[
b(p,u^2,u) = \max_{\substack{q\mid u\\q\neq p}} \left\{\nu_p(q^{p-1}-1) + \nu_p(\ord_{u_q}(q))\right\}.
\]
If $F(u^2,u)<u^2$ then $F(u^2,u)\leq u^2/p$ for some prime $p\mid u$, and $F(u^2,u)\geq u\varphi(u) = u^2\prod_{q\mid u}(1-1/q)$ implies that $\prod_{q\mid u}(1-1/q)^{-1} \geq p \geq 5$.
However, the product on the left here cannot exceed $1.7$ for $u\leq U$, so we require that $F(u^2,u)=u^2$.
Thus, for each prime $p\mid u$, there must exist another prime $q\mid u$ such that either $q^{p-1}\equiv1$ mod $p^2$, that is, $(q,p)$ is a Wieferich prime pair, or $p\mid\ord_{u_q}(q)$.
The latter condition requires that $p\mid(r-1)$ for some prime $r\mid u$ with $r\neq q$.

Our search for allowable lengths in the Barker sequence problem begins then with the construction of a large directed graph $D(U)$, whose vertices are a subset of the primes $p\equiv1$ mod $4$ with $p\leq P$.
We place a \textit{solid edge} $q\solid p$ if $(q,p)$ is a Wieferich prime pair, and a \textit{flimsy edge} $r\flimsy p$ if $p\mid(r-1)$.
Thus, a solid edge $q\solid p$ indicates that the presence of $q$ as a divisor of $u$ guarantees that $b(p,u^2,u)\geq2$, and a flimsy edge $r\flimsy p$ means that it is possible that $b(p,u^2,u)\geq2$ if $r\mid u$, as long as $u$ has at least one other prime divisor.
(In fact, it is quite likely that $b(p,u^2,u)\geq2$ in this case if $u$ has several other prime factors, or if $p$ is large.)
We describe the construction and analysis of our directed graph $D(U)$ in three stages.

\subsection{Graph construction}\label{subsecConstruct}
We begin with the directed graph $D(U_0)$ constructed in \cite{Moss09}, where $U_0=10^{15}/\sqrt{2}$, and then add the vertices and edges required for $D(U)$.
Most of the computations here were dedicated to detecting new Wieferich prime pairs $(q,p)$ where $q<p$.
By considering our new bound $U$ and prime bound $P$, the searches summarized in Table~\ref{tableAscWP} were required.
This search finds $156\,927$ new Wieferich pairs, involving $308\,837$ primes.
These pairs are available at the web site \cite{MossinghoffWeb}.
We remark that this computation required by far the largest portion of the CPU time required in this project, totaling approximately $7.5$ core-years.

\begin{table}[tbp]
\caption{Wieferich pair searches, $p\equiv q\equiv 1$ mod~$4$, $q<p$.}\label{tableAscWP}
\begin{center}
\begin{tabular}{cccc}
$q$ Range & $p$ Range & Pairs\\\hline
$[10^3,10^4]$\fixspacetop & $[10^{11}, P]$ & $3$\\
$[10^4,10^5]$ & $[10^{10}, P]$ & $234$\\
$[10^5,10^6]$ & $[7.5\cdot10^9, \min\{P,U/q\}]$ & $1537$\\
$[10^6,10^7]$ & $[7.5\cdot10^8, U/q]$ & $14981$\\
$[10^7,10^8]$ & $[10^8, U/q]$ & $125706$ &\\
$[10^8,\sqrt{U}]$ & $[q, \sqrt{U}]$ & $14466$\\\hline
\end{tabular}
\end{center}
\end{table}

After this, we find all Wieferich pairs $(q,p)$ with $q>p$, where $p\equiv1$ mod $4$ and $q$ is a prime that just appeared in this last search.
We then add a flimsy link for each prime $r\equiv1$ mod $4$ with $r\mid(q-1)$, for any newly appearing prime $q$,
This process continues, finding new descending Wieferich pairs and new flimsy links, until no additional new primes arise.
At the end of this process, our graph $D(U)$ has $608\,246$ vertices, $950\,456$ solid edges, and $665\,640$ flimsy edges.
For comparison, the graph $D(U_0)$ in \cite{Moss09} had $252\,905$ vertices, $387\,444$ solid edges, and $284\,272$ flimsy edges.

\subsection{Cycle enumeration}\label{subsecCycEnum}
Using Tarjan's algorithm \cite{Tarjan}, we find that there are $4656$ cycles in $D(U)$, with lengths from $2$ to $50$, and these are available at \cite{MossinghoffWeb}.
This is the same number that are present in the smaller graph $D(U_0)$.
We remark that only $2688$ cycles were reported for this graph in \cite{Moss09}, owing to an error in the implementation of Tarjan's algorithm.
However, none of the missing cycles was sufficiently short to alter the conclusions of \cite{Moss09}.

\subsection{Cycle augmentation and final processing}\label{subsecCycAug}
Only five cycles $C$ in $D(U)$ have the property that $\prod_{p\in C}p\leq U$: these are the same five cycles identified in \cite{Moss09} for $D(U_0)$.
We apply the same algorithm described in that paper to determine, for each such cycle $C$, all connected subgraphs $G$ of $D(U)$ containing $C$ for which $\prod_{p\in G}p\leq U$ and for which each $p\in G$ can be reached in $G$ from some vertex $q\in C$.
After this, we check that each such augmented cycle $G$ produces a viable value for $u$, by checking the required condition at each flimsy link.
We also check all of the values of the $b(p,u^2,u)$ in order to detect any admissible non-squarefree multiples of the allowed squarefree ones which lie below our bound $U$.

This produces sixteen permissible values for $u\leq U$, nine of which were reported in \cite{Moss09}.
There are seven new values:
\begin{align}
1\,087\,601\,914\,767\,745 &= 5 \cdot 13 \cdot 41 \cdot 2953 \cdot 138200401,\label{alnV1}\\
1\,258\,257\,961\,850\,525 &= 5^2 \cdot 193 \cdot 4877 \cdot 53471161,\label{alnV2}\\
2\,426\,188\,886\,789\,585 &= 5 \cdot 29 \cdot 41 \cdot 2953 \cdot 138200401,\label{alnV3}\\
5\,032\,969\,334\,448\,665 &= 5 \cdot 5333 \cdot 188748146801,\label{alnV4}\\
6\,308\,091\,105\,652\,921 &= 13 \cdot 29 \cdot 41 \cdot 2953 \cdot 138200401,\label{alnV5}\\
13\,337\,534\,395\,615\,565 &= 5 \cdot 53 \cdot 193 \cdot 4877 \cdot 53471161,\label{alnV6}\\
31\,540\,455\,528\,264\,605 &= 5 \cdot 13 \cdot 29 \cdot 41 \cdot 2953 \cdot 138200401.\label{alnV7}
\end{align}

We now turn to the remaining three criteria from Section~\ref{secRestrictions}.
The self-conjugacy test of Theorem~\ref{thmTuryn} rules out \eqref{alnV2}, \eqref{alnV4}, and \eqref{alnV6}.
(For \eqref{alnV2}, use $r=53471161$ and $s=2 \cdot 4877^2 r^2$, for \eqref{alnV4}, use $r=5333$ and $s=188748146801^2r^2$, and for \eqref{alnV6}, use $r=4877 \cdot 53471161$ and $s=r^2$.)
Theorem~\ref{thmLS5} then disqualifies \eqref{alnV1}, \eqref{alnV3}, and \eqref{alnV5}.
(Use $p=138200401$ in all three cases; for \eqref{alnV1}, use $r=5 \cdot 2953$, for \eqref{alnV3} and \eqref{alnV5}, use $r=29 \cdot 2953$.)
The last possibility \eqref{alnV7} cannot be excluded by using Theorem~\ref{thmTuryn}, Theorem~\ref{thmLS5}, or Theorem~\ref{thmLS1}.
This establishes Theorem~\ref{thmBarkerBound}.\qed

\section{Additional admissible values}\label{secAdditional}

We may use our graph $D(U)$ to construct some additional admissible values for the length of a Barker sequence by raising the threshold used in Section~\ref{subsecCycAug} from $U=10^{16.5}$ to some larger value, $W$.
Since we perform no additional searches for Wieferich prime pairs as in Section~\ref{subsecConstruct} and do not construct $D(W)$, we cannot conclude that we find all permissible lengths of Barker sequences in the range $[4U^2, 4W^2]$.
However, we describe these additional admissible values here in order to facilitate future research efforts.

\subsection{50 digits}\label{subsec50digits}
We set $W=5\cdot10^{24}$, and use the method of Section~\ref{subsecCycAug} to obtain a number of integers $u\in[U,W]$ that pass Corollary~\ref{corLargePrime}, Theorem~\ref{thmLS}, and Theorem~\ref{thmEKS}.
We find $133$ such integers, and we apply the other three criteria to these values.
Theorem~\ref{thmTuryn} disqualifies $115$ of them, and Theorems~\ref{thmLS5} and~\ref{thmLS1} do not exclude any of the remaining values.
Table~\ref{table50} summarizes the results of each of these latter three criteria indexed by $\Omega(u)$, which denotes the total number of prime factors of $u$, counting multiplicity.
This leaves eighteen admissible values for $u\in[U,W]$, and therefore eighteen additional plausible integers $n<10^{50}$ for the length of a long Barker sequence.
These values are listed in Table~\ref{tableAdmissible50}.

\begin{table}[tbp]
\caption{Eighteen admissible $u\in[10^{16.5}, 5\cdot10^{24}]$.}\label{tableAdmissible50}
\begin{center}
\begin{tabular}{r@{\qquad}l}
\multicolumn{1}{c}{$u$} & Factorization\\\hline
66\,687\,671\,978\,077\,825\fixspacetop & $5^2 \cdot 53 \cdot 193 \cdot 4877 \cdot 53471161$\\
866\,939\,735\,715\,011\,725 & $5^2 \cdot 13 \cdot 53 \cdot 193 \cdot 4877 \cdot 53471161$\\
1\,293\,740\,836\,374\,709\,805 & $5 \cdot 53 \cdot 97 \cdot 193 \cdot 4877 \cdot 53471161$\\
6\,468\,704\,181\,873\,549\,025 & $5^2 \cdot 53 \cdot 97 \cdot 193 \cdot 4877 \cdot 53471161$\\
16\,818\,630\,872\,871\,227\,465 & $5 \cdot 13 \cdot 53 \cdot 97 \cdot 193 \cdot 4877 \cdot 53471161$\\
84\,093\,154\,364\,356\,137\,325 & $5^2 \cdot 13 \cdot 53 \cdot 97 \cdot 193 \cdot 4877 \cdot 53471161$\\
2\,487\,505\,958\,525\,418\,181\,705 & $5 \cdot 29 \cdot 41 \cdot 2953 \cdot 1025273 \cdot 138200401$\\
6\,467\,515\,492\,166\,087\,272\,433 & $13 \cdot 29 \cdot 41 \cdot 2953 \cdot 1025273 \cdot 138200401$\\
19\,417\,213\,258\,149\,231\,605\,065 & $5 \cdot 17 \cdot 613 \cdot 1974353 \cdot 188748146801$\\
32\,337\,577\,460\,830\,436\,362\,165 & $5 \cdot 13 \cdot 29 \cdot 41 \cdot 2953 \cdot 1025273 \cdot 138200401$\\
863\,383\,081\,390\,130\,269\,759\,645 & $5 \cdot 41 \cdot 193 \cdot 2953 \cdot 53471161 \cdot 138200401$\\
1\,686\,504\,775\,565\,176\,744\,556\,405 & $5 \cdot 13 \cdot 29 \cdot 41 \cdot 2953 \cdot 53471161 \cdot 138200401$\\
1\,890\,448\,348\,089\,674\,770\,182\,781 & $53 \cdot 97 \cdot 4794006457 \cdot 76704103313$\\
2\,630\,496\,319\,975\,038\,327\,042\,325 & $5^2 \cdot 193 \cdot 24697 \cdot 53471161 \cdot 412835053$\\
2\,988\,996\,856\,098\,832\,119\,836\,165 & $5 \cdot 13 \cdot 123397 \cdot 1974353 \cdot 188748146801$\\
3\,080\,894\,677\,428\,239\,302\,747\,085 & $5 \cdot 5333 \cdot 612142549 \cdot 188748146801$\\
3\,770\,469\,237\,344\,599\,632\,723\,365 & $5 \cdot 53 \cdot 97 \cdot 193 \cdot 4877 \cdot 2914393 \cdot 53471161$\\
4\,316\,915\,406\,950\,651\,348\,798\,225 & $5^2 \cdot 41 \cdot 193 \cdot 2953 \cdot 53471161 \cdot 138200401$\\\hline
\end{tabular}
\end{center}
\end{table}

\begin{table}[tbp]
\caption{Effect of Theorems~\ref{thmTuryn}, \ref{thmLS5}, and~\ref{thmLS1} for $W=5\cdot10^{24}$.}\label{table50}
\begin{center}
\begin{tabular}{c|c|ccc|c}
&\multicolumn{1}{c}{Initial}&\multicolumn{3}{c}{Exclusions:}&\multicolumn{1}{c}{Admissible}\\
$\Omega(u)$ & Number & Thm~\ref{thmTuryn} & Thm~\ref{thmLS5} &Thm~\ref{thmLS1} & Number\\\hline
3\fixspacetop & 7 & 7 & -- & -- & 0\\
4 & 27 & 25 & 0 & 0 & 2\\
5 & 46 & 44 & 0 & 0 & 2\\
6 & 41 & 35 & 0 & 0 & 6\\
7 & 11 & 4 & 0 & 0 & 7\\
8 & 1 & 0 & 0 & 0 & 1\\\hline
Total\fixspacetop & 133 & 115 & 0 & 0 & 18
\end{tabular}
\end{center}
\end{table}

\subsection{100 digits}\label{subsec100digits}
We set $W=5\cdot10^{49}$, and in the same way we obtain $238\,858$ values of $u\leq W$ from the graph $D(U)$ which satisfy Corollary~\ref{corLargePrime}, Theorem~\ref{thmLS}, and Theorem~\ref{thmEKS}.
This includes the $133$ values from Section~\ref{subsec50digits}, as well as those from Section~\ref{secComputations} and from \cite{Moss09}.
However, because our tests for Theorems~\ref{thmTuryn} and~\ref{thmLS1} have running time that is exponential in $\Omega(u)$, only a portion of these could be examined against the remaining criteria.
In particular, we checked Theorem~\ref{thmTuryn} for $u$ with $\Omega(u)\leq8$, and Theorem~\ref{thmLS1} for $\Omega(u)\leq6$.
On the other hand, we tested Theorem~\ref{thmLS5} on all the values, since very few had a prime-power factor large enough to qualify in this theorem.
In all, we ruled out $1051$ of these values by using these last three criteria, nearly all coming from Turyn's test.
Table~\ref{table100} summarizes these results.
Combined with Table~\ref{table50}, we find $549$ integers $n<10^{100}$ that satisfy all of the known requirements for the length of a Barker sequence.

\begin{table}[tbp]
\caption{Effect of Theorems~\ref{thmTuryn}, \ref{thmLS5}, and~\ref{thmLS1} for $W=5\cdot10^{49}$.}\label{table100}
\begin{center}
\begin{tabular}{c|c|ccc|c}
&\multicolumn{1}{c}{Initial}&\multicolumn{3}{c}{Exclusions:}&\multicolumn{1}{c}{Admissible}\\
$\Omega(u)$ & Number & Thm~\ref{thmTuryn} & Thm~\ref{thmLS5} &Thm~\ref{thmLS1} & Number\\\hline
2\fixspacetop & 1 & 1 & -- & -- & 0\\
3 & 10 & 10 & -- & -- & 0\\
4 & 48 & 44 & 1 & 0 & 3\\
5 & 185 & 117 & 3 & 0 & 65\\
6 & 701 & 226 & 0 & 2 & 473\\
7 & 2560 & 326 & 0 && $\leq2234$\\
8 & 8440 & 321 & 0 && $\leq8119$\\
9 & 22406 && 0 && $\leq22406$\\
10 & 43523 && 0 && $\leq43523$\\
11 & 59673 && 0 && $\leq59673$\\
12 & 55200 && 0 && $\leq55200$\\
13 & 32627 && 0 && $\leq32627$\\
14 & 11266 && 0 && $\leq11266$\\
15 & 2029 && 0 && $\leq2029$\\
16 & 168 && 0 && $\leq168$\\
17 & 21 && 0 && $\leq21$\\\hline
Total\fixspacetop & 238858 & 1045 & 4 & 2 & $\leq237807$
\end{tabular}
\end{center}
\end{table}

In addition, we note that Turyn's self-conjugacy test is empirically less effective as $\Omega(u)$ grows larger.
Indeed, Theorem~\ref{thmTuryn} disqualifies $91.7\%$ of the values in Table~\ref{table100} with $\Omega(u)=4$, but only $3.8\%$ of those with $\Omega(u)=8$.
In the same way, we expect that Theorem~\ref{thmLS1} will grow less effective with larger $\Omega(u)$, owing to the embedded self-conjugacy test.
It seems reasonable to expect then that the vast majority of the remaining $237\,258$ integers $u$ from Table~\ref{table100} will pass these two criteria.
All $237\,807$ values are available at the authors' web site \cite{MossinghoffWeb} to assist future research efforts.

We note that the largest surviving integer $u<5\cdot10^{49}$ is
\[
u = 49\,998\,876\,926\,572\,332\,623\,608\,513\,080\,060\,406\,791\,494\,480\,564\,249,
\]
whose factorization is $13 \cdot 29 \cdot 37 \cdot 41 \cdot 53 \cdot 89 \cdot 97 \cdot 149 \cdot 349 \cdot 3049 \cdot 12197 \cdot 268693 \cdot 4794006457 \cdot 76704103313$, and that the largest such $u$ with $\Omega(u)=17$ is
\[
u = 36\,333\,506\,323\,649\,215\,674\,622\,109\,967\,025\,227\,563\,243\,374\,636\,565,
\]
with factorization $5 \cdot 13^2 \cdot 17 \cdot 29 \cdot 37 \cdot 53 \cdot 97 \cdot 197 \cdot 653 \cdot 1381 \cdot 1777 \cdot 4057 \cdot 11821 \cdot 16229 \cdot 24329 \cdot 76704103313$.

Also, it is interesting to note that each of the $237\,807$ different values of $u$ remaining in Table~\ref{table100} contains at least one of the following seven cycles in $D(U)$:
\begin{align*}
&5 \solid 53471161 \flimsy 5, \quad 5 \solid 188748146801 \solid 5, \quad 41 \solid 138200401 \flimsy 2953 \flimsy 41,\\
&53 \solid 97 \solid 76704103313 \flimsy 4794006457 \flimsy 53,\\
&30109 \solid 1128713 \solid 268813277 \flimsy 2167849 \flimsy 30109,\\
&37 \solid 76407520781 \flimsy 3301 \solid 24329 \solid 1297 \solid 31268910217 \solid 2797 \solid 76369 \flimsy 37,\\
&53 \solid 97 \solid 76704103313 \solid 16229 \flimsy 4057 \solid 11821 \flimsy 197 \solid 653 \solid 1381 \solid 1777 \solid 53.
\end{align*}
Many other cycles occur: some that contain one of these seven, such as $5 \solid 53471161 \solid 193 \solid 5$, and some that appear only in combination with one of them, like $5 \solid 6692367337 \solid 1601 \solid 5$.

\section{Circulant Hadamard matrices}\label{secCHM}

We conclude with a brief update on computational results in the circulant Hadamard matrix problem.
In \cite{Moss09}, the second author used the method of Section~\ref{secComputations}, but including the primes congruent to $3$ mod $4$, to find all integers $u\leq U_0'=10^{13}$ that satisfy Corollary~\ref{corLargePrime} and Theorem~\ref{thmLS}.
This effort produced a list of $2064$ integers.
However, this computation missed a number of cycles in the directed graph that was constructed there, due to the error in the earlier implementation of Tarjan's algorithm.
In \cite{Moss09}, $461\,653$ distinct cycles with length at most $12$ were reported; in fact, there are $1\,126\,465$ such cycles.
Only $196$ of these have vertex product less than $U_0'$ (up from $154$ in the earlier work), leading to $7676$ values after the augmenting step (up from $7491$).
However, after the final filtering steps, involving checking the bound on $F$ from Theorem~\ref{thmLS} and checking for valid non-squarefree values, we obtain precisely the same set of $2064$ integers satisfying Corollary~\ref{corLargePrime} and Theorem~\ref{thmLS}.

In addition, the same article reported that $486$ of these $2064$ values were eliminated by using Theorem~\ref{thmTuryn}.
More is true: in fact, $611$ of these values are disqualified by using Turyn's self-conjugacy test.
The software used to implement Turyn's test in \cite{Moss09} was unable to resolve a number of cases; a new implementation using Sage performed much better.

Finally, in their recent paper \cite{LS12}, Leung and Schmidt obtained a further restriction for the circulant Hadamard matrix problem, covering a special case that does not include the Barker sequence question.
We state this result here.

\begin{thm}\label{thmLS10}
Let $u$ be an integer whose prime divisors are all congruent to $3$ mod $4$.
Let $p$ be one of these divisors, and suppose that $w$ is a divisor of $u$ that is self-conjugate modulo $p$.
Let $q_1$, \ldots, $q_k$ be the prime divisors, excluding $p$, of $u/w$.
If $u=w$ or
\[
\gcd(\ord_p(q_1),\ldots,\ord_p(q_k)) \leq \frac{u^2}{w^2},
\]
then no circulant Hadamard matrix of order $4u^2$ exists.
\end{thm}

They used this result, together with Theorems~\ref{thmLS5} and~\ref{thmLS1}, to establish that $185$ of the remaining values from \cite{Moss09} are not allowed as the order of a circulant Hadamard matrix.
We find that $103$ of these $185$ values may also be eliminated by using Turyn's criterion.

Table~\ref{tableCHM} summarizes the results of applying Theorems~\ref{thmTuryn}, \ref{thmLS5}, \ref{thmLS1}, and \ref{thmLS10} successively to the $2064$ integers $u\leq U_0'$ that satisfy Corollary~\ref{corLargePrime} and Theorem~\ref{thmLS}.
Because Turyn's criterion now accounts for $22$ of the values not shown to be inadmissible in \cite{LS12} or \cite{Moss09}, the number of permissible orders $n$ for a circulant Hadamard matrix, with $4<n\leq4\cdot10^{26}$, now stands at $1371$.

\begin{table}[tbp]
\caption{Effect of Theorems~\ref{thmTuryn}, \ref{thmLS5}, \ref{thmLS1}, and~\ref{thmLS10} for $U_0'=10^{13}$.}\label{tableCHM}
\begin{center}
\begin{tabular}{c|c|cccc|c}
&\multicolumn{1}{c}{Initial}&\multicolumn{4}{c}{Exclusions:}&\multicolumn{1}{c}{Admissible}\\
$\Omega(u)$ & Number & Thm~\ref{thmTuryn} & Thm~\ref{thmLS5} &Thm~\ref{thmLS1} &Thm~\ref{thmLS10} & Number\\\hline
2\fixspacetop & 6 & 6 & -- & -- & -- & 0\\
3 & 50 & 42 & 0 &0 & 0 & 8\\
4 & 216 & 149 & 4 & 5 & 2 & 56\\
5 & 496 & 200 & 11 & 9 & 7 & 269\\
6 & 644 & 147 & 3 & 12 & 11 & 471\\
7 & 453 & 56 & 1 & 6 & 10 & 380\\
8 & 170 & 11 & 0 & 0 & 1 & 158\\
9 & 28 & 0 & 0 & 0 & 0 & 28\\
10 & 1 & 0 & 0 & 0 & 0 & 1\\\hline
Total\fixspacetop & 2064 & 611 & 19 & 32 & 31 & 1371
\end{tabular}
\end{center}
\end{table}

We remark that there are presently just five allowable integers $u$ in the circulant Hadamard matrix problem for which $n=4u^2\leq 10^{14}$:
\begin{gather*}
11\,715 = 3 \cdot 5 \cdot 11 \cdot 71, \quad 82\,005 = 3 \cdot 5 \cdot 7 \cdot 11 \cdot 71, \quad 550\,605 = 3 \cdot 5 \cdot 11 \cdot 47 \cdot 71,\\
3\,854\,235 = 3 \cdot 5 \cdot 7 \cdot 11 \cdot 47 \cdot 71, \quad 3\,877\,665 = 3 \cdot 5 \cdot 11 \cdot 71 \cdot 331,
\end{gather*}
and thirteen additional values of $u$ for which $n\leq10^{16}$:
$5\,418\,777$,
$8\,515\,221$,
$9\,031\,295$,
$9\,047\,885$,
$10\,975\,393$,
$12\,663\,915$,
$14\,192\,035$,
$16\,256\,331$,
$27\,093\,885$,
$27\,143\,655$,
$29\,549\,135$,
$32\,926\,179$, and
$42\,576\,105$.
The complete list of $1371$ integers $u$ is available at the web site \cite{MossinghoffWeb}.

\section*{Acknowledgements}

We thank the high performance computing consortium WestGrid, a subdivision of Compute Canada, as well as the Centre for Interdisciplinary Research in the Mathematical and Computational Sciences (IRMACS) at Simon Fraser University, for computational resources.


\begin{bibdiv}
\begin{biblist}

\bib{EK92}{article}{
      author={Eliahou, S.},
      author={Kervaire, M.},
       title={Barker sequences and difference sets},
        date={1992},
     journal={Enseign. Math. (2)},
      volume={38},
      number={3-4},
       pages={345\ndash 382},
        note={Corrigendum, Enseign. Math. (2) \textbf{40} (1994), no. 1--2,
  109--111},
      review={\MR{MR1189012 (93i:11018)}},
}

\bib{EKS}{article}{
      author={Eliahou, S.},
      author={Kervaire, M.},
      author={Saffari, B.},
       title={A new restriction on the lengths of {G}olay complementary
  sequences},
        date={1990},
     journal={J. Combin. Theory Ser. A},
      volume={55},
      number={1},
       pages={49\ndash 59},
      review={\MR{MR1070014 (91i:11020)}},
}

\bib{JL}{article}{
      author={Jedwab, J.},
      author={Lloyd, S.},
       title={A note on the nonexistence of {B}arker sequences},
        date={1992},
     journal={Des. Codes Cryptogr.},
      volume={2},
      number={1},
       pages={93\ndash 97},
      review={\MR{MR1157481 (93e:11032)}},
}

\bib{LS05}{article}{
      author={Leung, K.~H.},
      author={Schmidt, B.},
       title={The field descent method},
        date={2005},
     journal={Des. Codes Cryptogr.},
      volume={36},
      number={2},
       pages={171\ndash 188},
      review={\MR{MR2211106 (2007g:05023)}},
}

\bib{LS12}{article}{
      author={Leung, K.~H.},
      author={Schmidt, B.},
       title={New restrictions on possible orders of circulant {H}adamard
  matrices},
        date={2012},
        ISSN={0925-1022},
     journal={Des. Codes Cryptogr.},
      volume={64},
      number={1-2},
       pages={143\ndash 151},
         url={http://dx.doi.org/10.1007/s10623-011-9493-1},
      review={\MR{MR2914407}},
}

\bib{Moss09}{article}{
      author={Mossinghoff, M.~J.},
       title={Wieferich pairs and {B}arker sequences},
        date={2009},
        ISSN={0925-1022},
     journal={Des. Codes Cryptogr.},
      volume={53},
      number={3},
       pages={149\ndash 163},
         url={http://dx.doi.org/10.1007/s10623-009-9301-3},
      review={\MR{MR2545689 (2011c:11039)}},
}

\bib{MossinghoffWeb}{misc}{
      author={Mossinghoff, M.~J.},
       title={Wieferich prime pairs, {B}arker sequences, and circulant
  {H}adamard matrices},
        date={2013},
        note={http://www.cecm.sfu.ca/$\sim$mjm/WieferichBarker},
}

\bib{Schmidt99}{article}{
      author={Schmidt, B.},
       title={Cyclotomic integers and finite geometry},
        date={1999},
     journal={J. Amer. Math. Soc.},
      volume={12},
      number={4},
       pages={929\ndash 952},
      review={\MR{MR1671453 (2000a:05042)}},
}

\bib{Tarjan}{article}{
      author={Tarjan, R.},
       title={Enumeration of the elementary circuits of a directed graph},
        date={1973},
     journal={SIAM J. Comput.},
      volume={2},
       pages={211\ndash 216},
      review={\MR{MR0325448 (48 \#3795)}},
}

\bib{Turyn65}{article}{
      author={Turyn, R.},
       title={Character sums and difference sets},
        date={1965},
     journal={Pacific J. Math.},
      volume={15},
       pages={319\ndash 346},
      review={\MR{MR0179098 (31 \#3349)}},
}

\bib{TS}{article}{
      author={Turyn, R.},
      author={Storer, J.},
       title={On binary sequences},
        date={1961},
        ISSN={0002-9939},
     journal={Proc. Amer. Math. Soc.},
      volume={12},
       pages={394\ndash 399},
      review={\MR{MR0125026 (23 \#A2333)}},
}

\end{biblist}
\end{bibdiv}
\end{document}